\documentclass{article}
\usepackage{arxiv}

\usepackage{graphicx} 
\usepackage{color}
\usepackage{comment}
\usepackage{algorithm,algorithmic}
\usepackage{graphicx}
\usepackage{hyperref}
\usepackage{amsmath}
\usepackage{color}
\usepackage{subcaption}  
\usepackage{amsmath} 
\usepackage[table]{xcolor}
\usepackage{natbib}

\author{
Bartol Borozan\\
Faculty of Informatics and Data Science \\
University of Regensburg, Germany \\
School of Applied Mathematics and Informatics \\
Josip Juraj Strossmayer University of Osijek, Croatia \\
\And
Alice Descoeudres\\
Faculty of Informatics and Data Science\\ 
University of Regensburg\\ Regensburg, Germany
\And
Aneta Neumann\\
Optimisation and Logistics\\
School of Computer Science \\ and Information Technology\\
Adelaide University\\
Adelaide, Australia
\And
Frank Neumann\\
Optimisation and Logistics\\
School of Computer Science \\ and Information Technology\\
Adelaide University\\
Adelaide, Australia
\And
Stefan Canzar\\
Faculty of Informatics and Data Science\\ 
University of Regensburg\\ Regensburg, Germany
}

\title{Combinatorial CAR T-cell Circuit Design using Single- and Multi-Objective Local Search}

\begin{document}

\maketitle

\begin{abstract}
CAR T-cell circuit design enables targeting tumour cells while keeping healthy cells unaffected to a cancer treatment. Designing CAR T-cell circuits that maximize the number of tumour cells eliminated (called efficacy) while meeting the safety constraint of keeping a large number of healthy cells unaffected is a challenging task. We provide a new multistart local search approach for this problem and show that it outperforms the current SOTA approach on a single-cell cancer atlas. Furthermore, we introduce a multi-objective local search approach which provides insights into the trade-offs of efficacy and safety of the treatment. Our experiments also demonstrate that the multi-objective approach yields improved solutions even for the constrained single-objective optimisation problem.
\end{abstract}

\section{Introduction}
T cells are an essential part of the adaptive immune system, recognising antigens and killing infected cells. These capabilities are harnessed for cancer immunotherapy. The treatment of patients with haematological malignancies, such as B cell lymphoma and multiple myeloma, has been revolutionised by genetically engineered T cells expressing chimeric antigen receptors (CARs), which combine targeted antigen binding and T-cell activation domains~\cite{melenhorst2022,sheykhhasan2024,zhang2017,qu2022}. T cells expressing these CAR constructs exhibit targeted cytotoxicity, thereby allowing tumour-specific T-cell populations to be created \textit{ex vivo} from patient-derived, non-tumour-specific T cells. The efficacy of CAR T-cell therapy depends on the specificity and affinity of the CAR to the targeted antigen, and on the abundance and accessibility of this antigen within the group of malignant cells~\cite{vo2025,sterner2021}. Tumours can evade CAR T cells by antigen escape, wherein the cellular composition shifts to deplete populations exhibiting the target antigen. Conversely, high abundance of the target antigen on non-malignant cells can cause ``on-target, off-tumour'' toxicity, compromising the safety of the treatment. 

One approach to mitigating these limitations of single-antigen targeting is to construct multi-specific or logic-gated CAR T-cells~\cite{davies2016}. Such cells express multiple binding domains targeted to distinct antigens, leading to T-cell activation only if a Boolean logical condition is fulfilled~\cite{han2019,hamieh2023}. Antigen selection remains a central challenge in this multi-targeted approach. Thousands of surface antigens can potentially serve as targets. The search space grows combinatorially when considering AND, OR, and NOT-gated CAR T-cell designs.

Several computational approaches have been proposed for \textit{in silico} CAR T-cell circuit design. In some studies, AND and NOT gated two-gene combinations were identified from bulk RNA sequencing data using simple mutual exclusivity criteria~\cite{hu2021b,mackay2020}. Dannenfelser~\textit{et al.} screened surface antigen gene pairs and triplets expressed in AND and NOT gated patterns using a cluster separation-based ranking~\cite{dannenfelser2020}. The method MadHitter was developed to find the minimum number of targets needed to meet safety and efficacy thresholds based on OR-gated approaches~\cite{ahmadi2022}. Kwon~\textit{et al.} introduced a convolutional neural network (CNN)-based selection algorithm for two-gene logic gates, incorporating AND, OR, and NOT gates, based on single-cell data~\cite{kwon2023}. \textit{LogiCAR designer}, a recently proposed method that uses an iterative genetic algorithm, searches the space of logical circuits involving up to five genes~\cite{madan2025}. 

In this paper, we explore the use of iterative search algorithms for the gate-based design of CAR T-cells. 
The combinatorial nature of the design space makes iterative search algorithms such as local search~\cite{aarts1997local} and evolutionary algorithms~\cite{DBLP:series/ncs/EibenS15} particularly attractive.
The key problem that we tackle is to design a Boolean circuit that eliminates as many tumour cells from a given set $T$ while minimizing off-target effects on healthy cells from a given set $H$. Due to practical constraints on circuit design, we follow the example of \cite{madan2025} and consider circuits with at most $5$ input variables. However, the number of genes that can be assigned to these variables may comprise several thousand. Since circuits are restricted to at most five variables, the space of Boolean circuit topologies can be explored exhaustively. However, selecting the optimal gene-to-variable assignment for a fixed circuit design constitutes the challenging and central optimisation problem. 

We first introduce a single-objective local search algorithm that maximizes the number of eliminated tumour cells while satisfying a safety constraint requiring a large fraction of healthy cells to remain unaffected. We show that a multistart variant consistently outperforms the genetic algorithm proposed in \textit{LogiCAR designer}~\cite{madan2025}. 

The constrained single-objective formulation works with a fixed threshold on the fraction of healthy cells that are not affected by the treatment. Multi-objective optimisation~\cite{DBLP:books/daglib/0022057,deb2001multi} allows one to explore trade-offs in terms of the given objective and constraints~\cite{DBLP:journals/nc/NeumannW06} and this approach has been shown to be highly beneficial to also obtain better solutions for constrained single-objective problems~\cite{DBLP:journals/ec/FriedrichN15,DBLP:books/sp/ZhouYQ19,DBLP:conf/ECJ/NeumannN20}.
We formulate a bi-objective optimisation problem that simultaneously maximizes tumour cell elimination and the number of healthy cells preserved. The resulting Pareto front provides insights into the efficacy–safety trade-off. Furthermore, we show that the multi-objective approach achieves even better solutions than the ones obtained by the multistart single-objective local search approach.

The remainder of this paper is organized as follows. Section~\ref{sec2} formally defines the problem and presents our single-objective formulation. Section~\ref{sec4} describes our single-objective local search algorithm, while Section~\ref{sec7} introduces the multi-objective local search algorithm. The datasets used in our evaluation are described in Section~\ref{sec5}, and experimental results are presented in Section~\ref{sec6}. Finally, Section~\ref{sec8} concludes the paper.

\section{Problem Definition}
\label{sec2}
The aim of antigen selection for multi-targeted CAR T-cell design is to identify a set of genes whose combinatorial expression, under a specific Boolean logic, maximizes the separation between tumour and normal (healthy) cell populations. CAR T-cells engineered to target this  combination of antigens are activated by as many tumour cells as possible, thereby killing the recognized cells, while limiting off-target effects to normal tissue to less than 10\%~\cite{kwon2023, mackay2020}. We formalize this objective as an optimisation problem as follows. 

\subsection{Formal Problem Definition}

Let $T\in\{0,1\}^{m_T\times n}$ and $H\in\{0,1\}^{m_H\times n}$ denote two binary expression matrices representing tumour and healthy cells, respectively.
The rows correspond to cells and the columns to genes. We consider Boolean circuits with $N$ input variables, where $N\leq 5$. The circuit topology (i.e., the arrangement of Boolean gates) is assumed to be fixed. The optimisation variables are the assignments of  
genes to the circuit input variables. 

Formally, let $f(x)=f(x_1, \ldots, x_N)$ be the Boolean function induced by the 
fixed circuit topology, where each input variable $x_i$ is assigned to one of the $n$ surfaceome genes. 
Let $f(x)\rvert_c$ be the value of this function when $x_i$ are set according to the (binary) expression values of cell $c$ (given by the corresponding rows of matrices $T$ or $H$). 
The goal is to assign surfaceome genes to the $N$ input variables of the Boolean circuit so that the number of tumour cells $t$ for which $f(x)\rvert_t=1$ is maximized, subject to the constraint that at least an $\alpha$-fraction of healthy cells $h$ satisfy $f(x)\rvert_h=0$.

We say a tumour cell $t$ is recognized by the circuit iff $f(x)\rvert_t=1$ holds.
We therefore seek an assignment $x=(x_1, \ldots, x_N)$ maximizing
\begin{equation}
F(x,T):=\sum_{t \in T}f(x)\rvert_t
\end{equation}
subject to
\begin{equation}
G(x, H) := \sum_{h \in H}(1-f(x)\rvert_h) \geq \alpha\cdot m_H
\end{equation}
holds.

Following previous studies, we use $\alpha=0.9$ and restrict the number of circuit inputs to $N \leq 5$.
In our datasets, the number of genes and cells are large, e.g. $n=3,000$, $m_T=27,000$, and $m_H=15,000$.

Given a solution $x=(x_1, \ldots, x_N)$, we require a fitness function $q(x, T, H)$ that quantifies its quality with respect to tumour and healthy cells and allows it to be compared with other candidate solutions.

\subsection{Data-based Fitness Evaluation}
Let $C(x,H) = \min\{G(x,H)- \alpha \cdot m_H,0\}$ measure the amount of constraint violation.
The term $C(x,H)$ is $0$ iff $x$ meets the constraint requirement and negative otherwise.
Among feasible solutions (meeting the constraint requirement), the goal is to maximize $F(x, T)$.

We define the fitness function $q(x, T, H) = (C(x,H), F(x,T))$ which is maximized in lexicographic order.
Hence, we have $q(y, T, H) \geq q(x, T, H)$ (in lexicographical order) iff
\begin{eqnarray}
& &(C(y, H) > C(x, H))\\ &\vee & 
((C(x, H) = C(y, H)) \wedge (F(y, T) \geq F(x, T)))
\end{eqnarray}
We use the fitness function $q(x,T,H)$ to direct the search.
We may write $q(x)$ instead of $q(x, T, H)$ if $T$ and $H$ are clear from the context.

\begin{algorithm}[t] 
\begin{algorithmic}

\REQUIRE{Set of candidate genes $X=\{x_1, \ldots, x_n\}$, binary matrices $T$ and $H$, definition of neighbourhood $Z(x)$}
\STATE Choose initial solution $x=(x_1, \ldots, x_N)$, where $x_i \in X$, $1 \leq i \leq N$.
\STATE Evaluate $q(x, T, H)$.
\STATE $t:=1$;
\WHILE{ $t \leq t_{\max}$}{
\STATE imp $=$ false;
\STATE $Z=Z(x)$;
\WHILE{$Z \not= \emptyset \wedge (!\mbox{imp})$}{
\STATE Choose solution $y \in Z$;
\STATE Evaluate $q(y, T, H)$.
\STATE $Z:= Z \setminus \{y\}$
\IF{$q(y, T, H) > q(x, T, H)$}{
\STATE $x:=y$, imp $=$ true}
\ENDIF
}
\ENDWHILE
\STATE $t:=t+1$;
}
\ENDWHILE
\STATE $x^*:=x$;
\STATE Output $x^*$
\end{algorithmic}

\caption{Local Search (LS)}
\label{alg:ls}
\end{algorithm}

\section{Multistart Single-Objective Local Search}
\label{sec4}

Given the description of $f$ through the circuit gates, the goal is to find a feasible $x$, which maximizes $F(x,T)$. Our single-objective local search (LS) algorithm (see Algorithm~\ref{alg:ls}) starts with an initial assignment of the $N$ variables. In each iteration, a new solution $y$ is produced from the neighbourhood $Z(x)$ of the current solution $x$.
We consider the simple swap neighbourhood, in which one variable in the current solution is replaced by a variable not already contained in the solution. Formally, for a solution $x=(x_1, \ldots, x_N)$,
the neighbourhood $Z(x)$ consists of all solutions obtained by replacing any input variable $x_i$, $1 \leq i \leq N$, with any gene $g\in X\setminus \{x_1, \ldots, x_n\}$, where $X$ denotes the set of all candidate genes. 
A neighbouring solution $y\in Z(x)$ is accepted if its fitness is strictly greater than that of the currentl solution $x$. The search terminates after exploring the neighbourhood of at most $t_{\max}$ solutions, and the best solution found is returned.

Since the quality of the final solution obtained by LS may depend strongly on the initial solution, we employ a multistart strategy. 
For each restart, the $N$ variables of the initial solution $x=(x_1, \ldots, x_N)$ are selected uniformly at random, without replacement, from the set of $n$ candidate variables. Unless stated otherwise, the multistart LS performs 100 independent runs, each starting from a newly generated random initial solution, and returns the best solution found.

\subsection{Neighbourhood and GPU-based Evaluation}

Since matrices $H$ and $T$ are generally large, efficient evaluation of the fitness function  $q(x, T, H)$ is crucial. Furthermore, local search benefits from evaluating multiple neighbouring solutions in parallel. After transferring the data to the GPU, we therefore evaluate neighbouring solutions in batches, allowing many candidate solutions to be processed simultaneously and substantially reducing the overall runtime.

\section{Multi-Objective Local Search}
\label{sec7}
We now develop a multi-objective formulation for the problem.
Rather than searching for a single solution satisfying a fixed safety threshold, we simultaneously optimize efficacy and safety. This enables us to explore the trade-off between these objectives, revealing how much efficacy can be gained by relaxing the safety requirement, or conversely, how much efficacy is sacrificed when requiring a higher level of safety.
Beyond providing insight into these trade-offs, the multi-objective formulation also induces a more flexible search process. Previous studies have shown that reformulating constrained single-objective problems as multi-objective optimization problems often leads to improved solutions for the original constrained problem.

Let $S(x)$ be the safety and $F(x)$ be the fitness of solution $x$. We define the bi-objective function $G(x) = (F(x), S(x))$ where we aim to maximize both objectives.
We say that $x$ weakly dominates $y$ (denoted as $x \succeq y$) iff $F(x) \geq F(y)$ and $S(x) \geq S(y)$.
We say that $x$ strongly dominates $y$ (denoted as $x \succ y$) iff $x \succeq y$ and $(F(x) \not = F(y) \vee S(x) \not = S(y))$.
Consistent with the single-objective formulation, we use 
$S(x) = G(x, H)$
and $F(x) = F(x,T)$.

\begin{algorithm}[t]
\begin{algorithmic}

\REQUIRE{Set of candidate genes $X=\{x_1, \ldots, x_n\}$, binary matrices $T$ and $H$, definition of neighbourhood $Z(x)$}
\STATE Choose initial solution $x=(x_1, \ldots, x_N)$, where $x_i \in X$, $1 \leq i \leq N$.
\STATE $A = \{x\}$\;
\STATE Evaluate $G(x)=(F(x), S(x))$.
\STATE $t:=1$;
\WHILE{ $t \leq t_{\max}$}{
\STATE imp $=$ false;
\STATE Choose $x \in A$ uniformly at random\;
\STATE $Z=Z(x)$;
\WHILE{$Z \not= \emptyset \wedge (!\mbox{imp})$}{
\STATE Choose solution $y \in Z$;
\STATE Evaluate $G(y)=(F(y), S(y))$.
\STATE $Z:= Z \setminus \{y\}$
\IF{there is no solution $z \in A$ with $z \succeq y$ }{
\STATE $A:= (A \setminus \{p \in A \mid y \succ p\}) \cup \{y\}$} 
(remove all solutions from $A$ strongly dominated by $y$ and add $y$ to $A$)
\ENDIF
\IF{$y \succ x$}{
\STATE imp $=$ true}
\ENDIF
}
\ENDWHILE
\STATE $t:=t+1$;
}
\ENDWHILE
\STATE $x^* = \arg \max_{x \in A} \{F(x) \mid S(x) \geq \alpha \cdot m_H\}$
\STATE Output $(A, x^*)$
\end{algorithmic}
\caption{Multi-Objective Local Search (MOLS)}
\label{alg:mols}
\end{algorithm}

Multi-Objective Local Search (MOLS) (see Algorithm~\ref{alg:mols}) extends our single-objective LS procedure to the multi-objective setting. The key difference is the use of an archive $A$ which stores one representative solution for each non-dominated objective vector discovered during the search.  
Since our goal is to characterize the complete trade-off between the two objectives, the archive is kept unbounded.

In each iteration, the algorithm selects a solution $x \in A$, explores its neighbourhood, and updates the archive accordingly. Given a newly generated solution $y$, it is inserted into $A$ iff no solution $z \in A$ weakly dominates it. If $y$ is added, all solutions in the archive that are strongly dominated by $y$ are removed.

In contrast to the single-objective setting, the search may switch between different regions of the search space by selecting different archive members for neighbourhood exploration. Whenever the neighbourhood of $x$ contains a solution $y$ that strongly dominates $x$, the current local search is terminated after updating the archive. A new solution is then selected from the archive, and its neighbourhood is explored.

The algorithm terminates after the neighbourhoods of at most $t_{\max}$ solutions have been explored. It returns the final archive $A$, representing the approximated Pareto front, together with the solution 
\begin{equation}
x^* = \arg\max_{x \in A,\; S(x)\ge \alpha m_H} F(x).
\end{equation}
i.e., the feasible solution with the highest efficacy according to the single-objective formulation.

\section{Experimental Setup and Datasets}
\label{sec5}
To evaluate the methods, we selected six major tumour types from the single-cell expression atlas built in~\cite{kwon2023} from 412 tumours and 12 normal organs: ovarian cancer (Ovary), non-small cell lung cancer (Lung), colorectal cancer (Colon), breast cancer (Breast), liver hepatocellular carcinoma (Liver), and pancreatic adenocarcinoma (Pancreas) (Table \ref{tab:dataset_summary}).
\begin{table}[t]
\setlength{\tabcolsep}{7pt}
\centering
\caption{Summary of datasets with cell counts and cancer cell proportions.}
\vspace{1mm}
\begin{tabular}{lcccc}
\hline
Dataset & \#Cells & \#Cancer & \#Healthy & \%Cancer \\
\hline
Breast    & 606720 & 517101 &  89619 & 85\% \\
Colon     & 455155 & 253475 & 201680 & 56\% \\
Liver     & 179067 & 138295 &  40772 & 77\% \\
Lung      & 488617 & 317505 & 171112 & 65\% \\
Ovary     & 248546 & 195539 &  53007 & 79\% \\
Pancreas  & 162932 & 111733 &  51199 & 69\% \\
\hline
\end{tabular}
\label{tab:dataset_summary}
\end{table}
Based on gene expression values, the datasets were converted into a binary format. All non-zero expression values were set to one, while zero values were retained as zero. 

We followed the strategy in~\cite{han2019} to evaluate and compare several methods to derive boolean formulas of sizes 2, 3, 4, and 5. We computed the safety and efficacy scores of the resulting formulas and compared their performance. The methods considered include LS (\ref{alg:ls}) in both regular and multistart variants, MOLS (\ref{alg:mols}), and a naive exhaustive search approach, in which all possible formulas are evaluated and the one with the highest efficacy score is selected subject to a safety score exceeding a predefined threshold. Additionally, we included a reimplementation of the state-of-the-art method \textit{LogiCAR designer}~\cite{madan2025}. \textit{LogiCAR designer} extended previous experimental work which considered up to 3 genes~\cite{han2019} and explored up to $N=5$ genes. Compared to previous computational approaches \cite{kwon2023} and \cite{dannenfelser2020} (limited to $N\leq 3$), it provided circuits with higher efficacy for a safety threshold of $90\%$. However, no public implementation of \textit{LogiCAR designer} was available. 
We therefore used our own implementation that follows the description provided in \cite{madan2025} as closely as possible, which we refer to as logicar$^*$. 

All methods were applied to a refined set of cancer surfaceome antigens, as described in \cite{kwon2023}. Following \cite{madan2025}, to reduce the number of possible circuit gates, we modified the datasets by adding, for each gene $g$, its negated counterpart $!g$. This eliminates the need for explicit NOT operations, thereby reducing the overall search space.
Additionally, we performed a 70\%/30\% train–test split. The methods were used to compute a boolean function on the 70\% training set, which was then evaluated on the remaining 30\% of the data. To ensure a fair comparison, all methods used the same training and testing split.

Each method was run 30 times to account for variability across runs. Exhaustive Search, being deterministic, was run only once. For each stochastic method, we report the mean efficacy and standard deviations across the 30 runs.
All experiments were run on an AMD Ryzen Threadripper 3990X 64-Core Processor @ 4.3GHz, 256 GB of RAM, an NVIDIA GeForce RTX 3090 GPU, and Python version 3.11.14. 

\section{Experimental Results}
\label{sec6}
We now investigate the performance of our approaches. We start by considering the constrained single-objective problem of maximizing efficacy. Afterwards, we provide multi-objective results showing trade-offs between efficacy and safety.

\subsection{Simple Circuits ($N=2,3$)}

We start by presenting experimental results for simple circuits.
Tables (\ref{tab:scores_2}) and (\ref{tab:scores_3}) show the results of all methods when searching for formulas of size 2 or 3. The tables report both the mean values across multiple runs and the standard deviation. The safety scores are not reported, as all methods are constrained to select only formulas with safety at or above 90\%, and in all experiments the resulting safety values remained within a 1\% margin of this threshold.
Due to the lower complexity of these cases, we also ran the exhaustive search algorithm, which provides the optimal efficacy score. 

\begin{table}[t]
\centering
\setlength{\tabcolsep}{5pt}
\begin{tabular}{l|cccc}
Dataset & \textbf{Local Search} & \textbf{logicar$^*$} & \textbf{Multistart LS} & \textbf{Exhaustive Search}\\
\hline
Breast & $65.42 \pm 5.1$ & $70.82 \pm 1.1$ & $71.22 \pm 0.0$ & $71.22$\\
Colon & $41.94 \pm 0.7$ & $42.40 \pm 0.9$ & $46.09 \pm 0.5$ & $46.21$\\
Liver & $41.84 \pm 1.3$ & $42.53 \pm 0.6$ & $42.90 \pm 0.0$ & $42.97$\\
Lung & $27.17 \pm 3.2$ & $29.19 \pm 0.8$ & $29.69 \pm 0.0$ & $29.69$\\
Ovary & $43.00 \pm 10.9$ & $54.24 \pm 1.1$ & $54.65 \pm 0.0$ & $54.65$\\
Pancreas & $64.06 \pm 5.6$ & $69.78 \pm 0.7$ & $70.71 \pm 1.1$ & $73.34$\\
\end{tabular}
\caption{Comparison of the efficacy of different methods for formula size 2, reported as mean $\pm$ standard deviation across runs. MOLS consistently achieves the optimal efficacy and is therefore omitted from the table.}
\label{tab:scores_2}
\end{table}

\begin{table}[t]
\centering
\setlength{\tabcolsep}{5pt}
\begin{tabular}{l|cccc}
Dataset & \textbf{Local Search} & \textbf{logicar$^*$} & \textbf{Multistart LS} & \textbf{Exhaustive Search}\\
\hline
Breast & $69.70 \pm 6.2$ & $77.02 \pm 0.0$ & $77.02 \pm 0.0$ & $77.02$\\
Colon & $44.32 \pm 3.2$ & $49.72 \pm 2.3$ & $50.58 \pm 0.0$ & $50.58$\\
Liver & $44.30 \pm 3.3$ & $48.41 \pm 0.0$ & $48.41 \pm 0.0$ & $48.41$\\
Lung & $29.12 \pm 1.8$ & $31.48 \pm 0.4$ & $32.08 \pm 0.0$ & $32.08$\\
Ovary & $52.28 \pm 10.5$ & $62.99 \pm 1.7$ & $64.66 \pm 0.0$ & $64.66$\\
Pancreas & $68.56 \pm 8.0$ & $80.52 \pm 1.2$ & $81.77 \pm 0.2$ & $82.72$\\
\end{tabular}
\caption{Comparison of the efficacy of different methods for formula size 3, reported as mean $\pm$ standard deviation across runs. MOLS consistently achieves the optimal efficacy and is therefore omitted from the table.}
\label{tab:scores_3}
\end{table}

The basic LS method exhibits the highest standard deviation due to the sensitivity to the choice of initial solution. The Multistart LS method, which performs multiple runs from different initial solutions, reduces the standard deviation and improves the mean performance. Both logicar$^*$ and Multistart LS achieve results close to, and in several cases equal to, the upper bound provided by exhaustive search. However, logicar$^*$ shows a slightly higher standard deviation across runs, whereas Multistart LS exhibits minimal or zero standard deviation in most cases. Notably, for some datasets (Pancreas and Colon), neither method consistently reaches the exhaustive search optimum. Additionally, we applied MOLS and report the solution with the highest efficacy among those satisfying the 90\% safety threshold. MOLS consistently attains the exhaustive search optimum and is therefore omitted from Tables (\ref{tab:scores_2}) and (\ref{tab:scores_3}).

\subsection{Complex Circuits ($N=4,5$)}

We now turn to more complex circuits, which are more challenging from an algorithmic perspective.
We performed the same experiments for formulas of sizes 4 and 5. However, due to the increased computational complexity, running the exhaustive search algorithm was no longer feasible. Tables \ref{tab:scores_4} and \ref{tab:scores_5} summarize the results, reporting the mean and standard deviation across runs. 

\begin{table}[t]
\centering
\setlength{\tabcolsep}{5pt}
\begin{tabular}{l|cccc}
Dataset & \textbf{Local Search} & \textbf{logicar$^*$} & \textbf{Multistart LS} & \textbf{Multi-Objective LS}\\
\hline
Breast & $75.33 \pm 3.8$ & $79.74 \pm 0.8$ & $80.04 \pm 0.0$ & $80.23 \pm 0.0$\\
Colon & $48.10 \pm 3.3$ & $51.34 \pm 0.9$ & $53.09 \pm 1.2$ & $55.83 \pm 0.0$\\
Liver & $47.42 \pm 3.3$ & $52.30 \pm 0.3$ & $52.50 \pm 0.1$ & $52.74 \pm 0.0$\\
Lung & $30.91 \pm 1.6$ & $32.65 \pm 0.3$ & $33.26 \pm 0.4$ & $35.42 \pm 0.0$\\
Ovary & $61.26 \pm 5.1$ & $67.55 \pm 0.9$ & $69.18 \pm 0.0$ & $68.92 \pm 0.6$\\
Pancreas & $77.12 \pm 5.3$ & $83.55 \pm 0.6$ & $84.34 \pm 0.4$ & $88.80 \pm 0.0$\\
\end{tabular}
\caption{Comparison of the efficacy of different methods for formula size 4, reported as mean $\pm$ standard deviation across runs.}
\label{tab:scores_4}
\end{table}

\begin{table}[t]
\centering
\setlength{\tabcolsep}{5pt}
\begin{tabular}{l|cccc}
Dataset & \textbf{Local Search} & \textbf{logicar$^*$} & \textbf{Multistart LS} & \textbf{Multi-Objective LS}\\
\hline
Breast & $78.40 \pm 1.5$ & $80.90 \pm 0.4$ & $81.24 \pm 0.2$ & $82.34 \pm 0.0$\\
Colon & $50.87 \pm 1.2$ & $53.15 \pm 0.5$ & $54.82 \pm 1.2$ & $58.08 \pm 0.1$\\
Liver & $50.05 \pm 2.7$ & $54.87 \pm 1.0$ & $55.68 \pm 0.0$ & $55.68 \pm 0.0$\\
Lung & $32.15 \pm 0.6$ & $33.05 \pm 0.2$ & $34.65 \pm 0.6$ & $36.85 \pm 0.0$\\
Ovary & $66.98 \pm 2.4$ & $70.81 \pm 0.5$ & $71.86 \pm 0.2$ & $72.48 \pm 0.0$\\
Pancreas & $82.33 \pm 1.6$ & $85.20 \pm 0.6$ & $87.21 \pm 1.1$ & $91.29 \pm 0.0$\\
\end{tabular}
\caption{Comparison of the efficacy of different methods for formula size 5, reported as mean $\pm$ standard deviation across runs.}
\label{tab:scores_5}
\end{table}

For formulas of sizes 4 and 5, MOLS achieves the highest mean efficacy across all datasets, outperforming both Multistart LS and logicar$^*$. A more detailed discussion on the results of the multi-objective approach is provided in Section~\ref{sec8}. Compared to MOLS and logicar$^*$, the basic LS method attains lower mean efficacy and exhibits greater variability across runs, although this variability decreases slightly as the formula size increases. In contrast, both Multistart LS and logicar$^*$ demonstrate low variability across runs, with Multistart LS often exhibiting near-zero standard deviation. The performance gap between Multistart LS and logicar$^*$ is generally small but consistent, with Multistart LS achieving a slightly higher mean efficacy than logicar$^*$ on every dataset.

\subsection{Trade-offs of Efficacy and Safety}
\label{sec8}
Instead of applying a fixed safety threshold of $90\%$, we now explore the trade-offs between efficacy and safety of the circuits designed by MOLS. 
Tables~\ref{tab:scores_4} and~\ref{tab:scores_5} show, for each dataset, the solution obtained by MOLS with the highest efficacy while satisfying the 90\% safety requirement. Consistent with previous findings, the search process underlying the multi-objective formulation frequently identifies better feasible solutions than the constrained single-objective formulation. For $N=2$ and $N=3$, selecting a formula based on the multi-objective search resulted in optimal solutions for all datasets. 
In contrast to~\cite{madan2025}, who observed only minor improvements in efficacy when increasing the number of genes beyond three, our multi-objective approach achieved substantial gains on several datasets. For example, efficacy increased from $82.72\%$ for $N=3$ to $91.29\%$ for $N=5$ on Pancreas, and from $50.58\%$ to $58.08\%$ on Colon.

\begin{figure}[t]
    \centering
    
    \includegraphics[width=.7\textwidth]{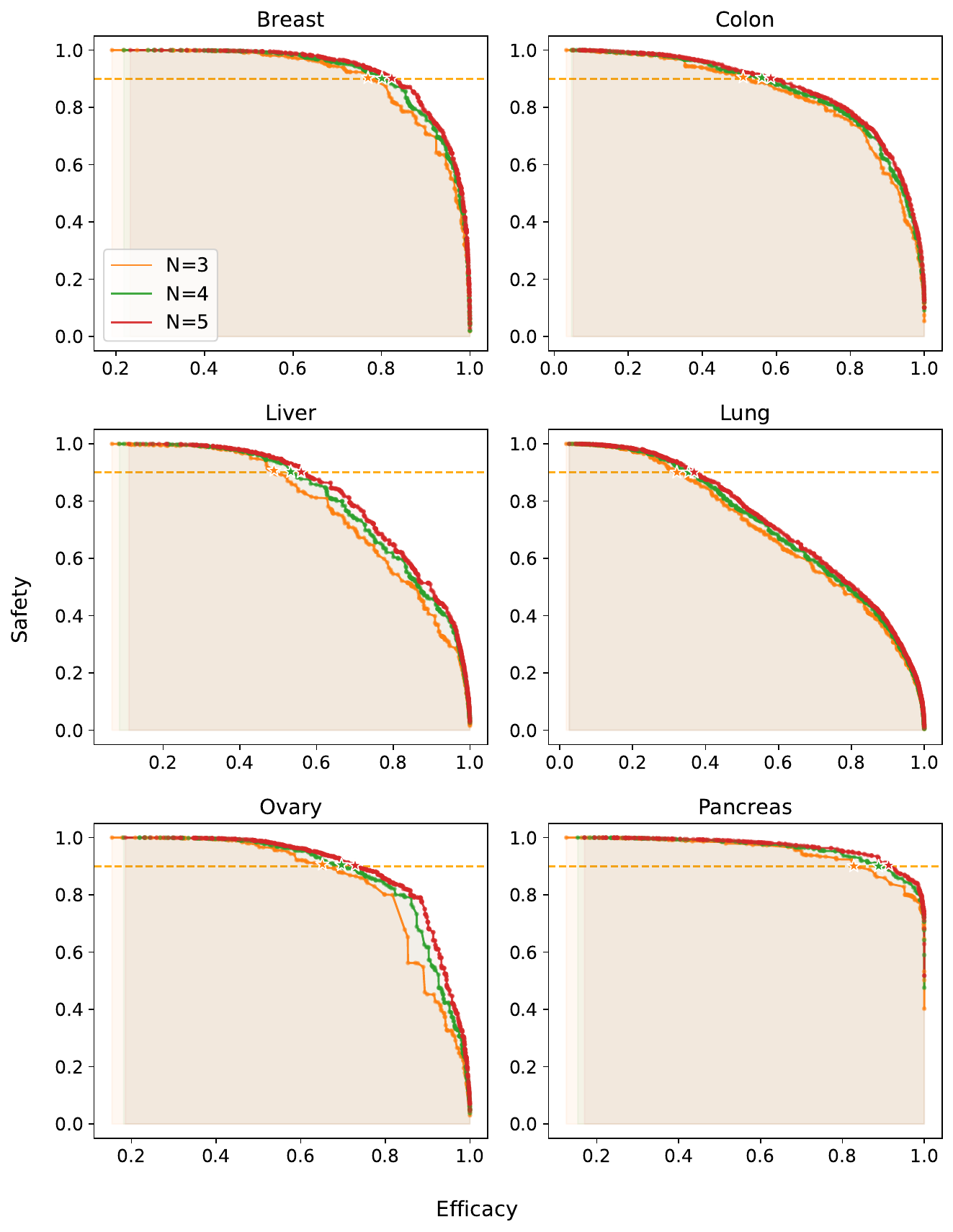}

    \caption{Pareto fronts obtained by MOLS for formulas of sizes $N=3, 4,$ and $5$. Each point represents a non-dominated solution in terms of efficacy and safety. The horizontal dashed line indicates the 90\% safety threshold. For each value of $N$, the feasible solution with the highest efficacy is highlighted.}
    \label{fig:multiobj_ls}
\end{figure}

Figure~\ref{fig:multiobj_ls} further illustrates the efficacy-safety trade-offs obtained by MOLS for formulas of sizes  $N=3, 4,$ and $5$. For several datasets, including Colon, Ovary, and Pancreas, the 90\% safety threshold lies just before the knee of the Pareto front, suggesting that modest reductions in safety could yield a substantial increase in efficacy. 
As the number of genes increases, the trade-off curves shift slightly towards higher efficacy, further supporting the observation that larger formulas provide better efficacy at a fixed safety threshold.

\subsection{Stability}
\begin{figure}
    \centering
    \includegraphics[width=1\linewidth]{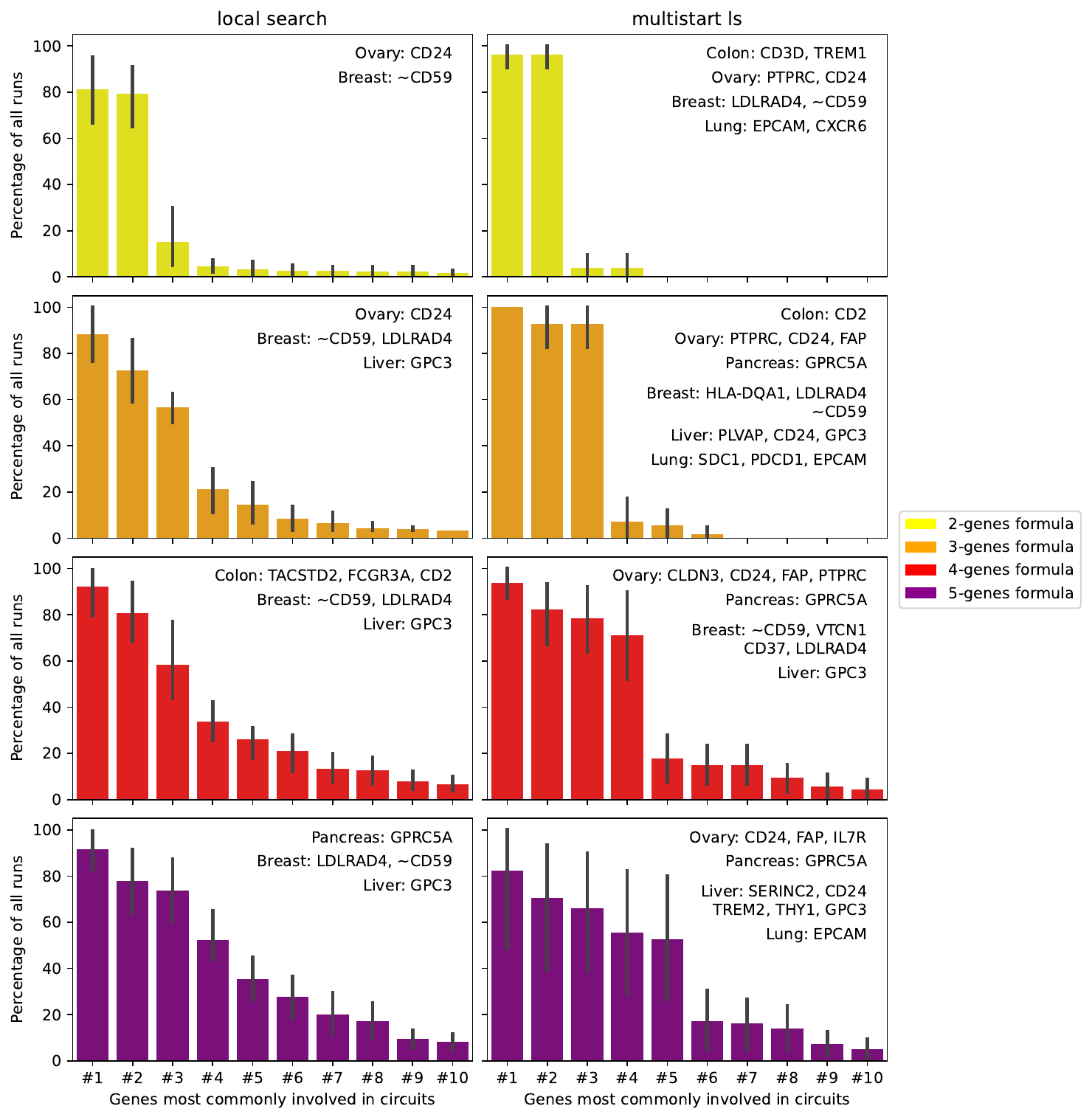}
    \caption{Stability of gene selection across formula sizes and tissues. Averaged across tissues, the frequency with which the 10 most frequently selected genes appear across all runs is shown for formulas of sizes $N=2$ (top) to $N=5$ (bottom). Core genes selected in 100\% of runs are highlighted for each tissue. A gene name preceded by $\sim$ indicates that the gene is included in every formula in negated form. Total number of runs: 30.}
    \label{fig:stability}
\end{figure}

Figure~\ref{fig:stability} illustrates the stability of LS and Multistart LS across datasets by measuring how consistently genes appear in highest-scoring formulas.

For LS, each of the two most frequently selected genes appears in approximately 80\% of two-gene formulas across all datasets. Multistart LS further increases this stability to nearly 100\%. Specifically, for the Colon, Ovary, Breast, and Lung datasets, every two-gene formula consists of the same pair of genes, with variation observed only for the Pancreas and Liver datasets. 

Similar trends are observed for larger formulas. For $N=3$ and $N=4$, Multistart LS identifies highly consistent gene sets, with the three and four most frequently selected genes appearing in more than 90\% and 75\% of formulas, respectively. For $N=5$, the stability of the core genes remains high for both methods, although additional genes vary more substantially across solutions. 

As expected, Multistart LS consistently produces more stable solutions than LS. For example, the immune regulator CD24 appears in every formula on the Ovary dataset generated by Multistart LS. While LS also identifies CD24 in all formulas for $N=2$ and $N=3$, it becomes less consistent for larger formulas, where more variable gene combinations are selected.

\subsection{Running Time}

\begin{figure}[t]
    \centering
    \makebox[\textwidth][c]{%
        \includegraphics[width=.65\textwidth]
         {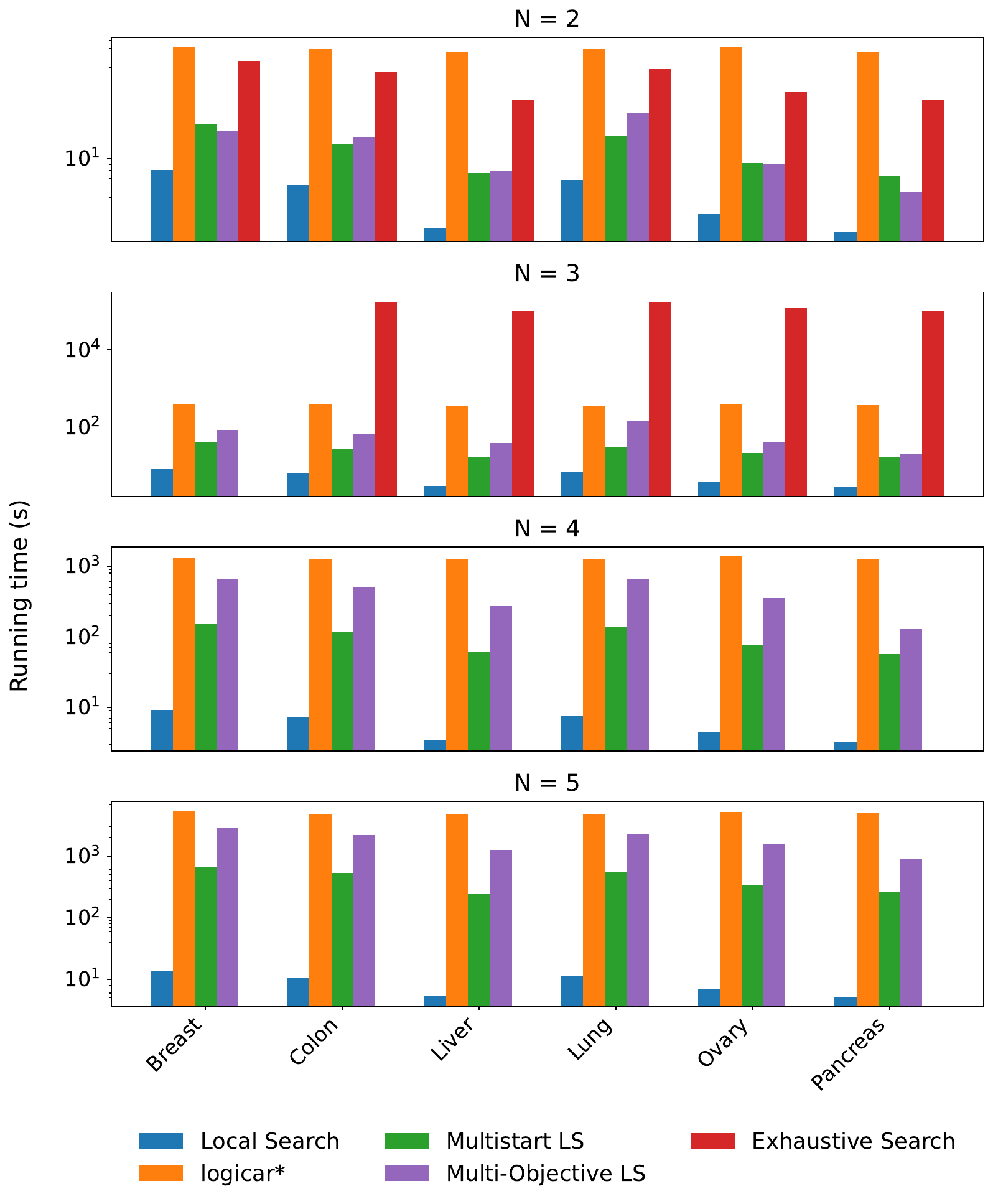}%
    }
    \caption{Average runtime (in seconds) of all methods evaluated in our experiments for formulas of sizes $N=2$ (top) to $N=5$ (bottom).}
    \label{fig:running_time}
\end{figure}

Figure~\ref{fig:running_time} compares the average runtime of the different methods across datasets as a function of formula size. 

LS is consistently the fastest method across all datasets, typically outperforming logicar$^*$ by two orders of magnitude or more (e.g., 14 seconds versus 4964 seconds on the Breast dataset for $N=5$). Multistart LS incurs additional runtime due to repeated restarts, but the overhead grows sublinearly with the number of restarts as CPU-side preprocessing and data preparation are performed only once and reused across all runs. In contrast, the runtime of MOLS approaches that of logicar$^*$ (2567 seconds versus 4965 seconds on the Breast dataset for $N=5$). As expected, exhaustive search is the slowest method for formulas of size $N=3$. However, for formulas of size $N=2$, its lower computational complexity, combined with efficient GPU utilisation, allows it to outperform logicar$^*$.

\section{Conclusions}

Combinatorial CAR T-cell circuit design is important for the treatment of various cancers. We have provided a simple Multistart LS approach that outperforms the recently proposed logicar$^*$ approach. Furthermore, we have developed a multi-objective model that allows to compute solution trade-offs with respect to efficacy and safety. The developed MOLS algorithm not only provides these additional insights but also improves upon the results that we obtained for the constrained single-objective formulation. While engineering CAR circuits targeting multiple genes remains experimentally challenging, our approach demonstrates that increasing the number of targets beyond three can yield further improvements in efficacy.

\bibliographystyle{abbrv}
\bibliography{ref} 

\end{document}